\documentclass[a4paper, 12pt] {article}
\usepackage[cp1251]{inputenc}
\usepackage[english,russian]{babel}
\usepackage{amsfonts}
\usepackage{mathtext}
\begin{document}

\thispagestyle{empty}
\begin{center}
{\bf SOLUTION OF THE INVERSE SPECTRAL PROBLEM\\ FOR A CONVOLUTION INTEGRO-DIFFERENTIAL OPERATOR\\ WITH ROBIN BOUNDARY CONDITIONS}
\end{center}

\begin{center}
{\bf S.A.\,Buterin\footnote{Department of Mathematics, Saratov University,
Astrakhanskaya 83, Saratov 410012, Russia, {\it Email: buterinsa@info.sgu.ru}} and
A.E.\,Choque Rivero\footnote{Instituto de F\'isica y Matem\'aticas, Universidad
Michoacana de San Nicol\'as de Hidalgo. Edificio C3A, Cd. Universitaria. C. P. 58040
Morelia, Mich., M\'exico, {\it Email: abdon@ifm.umich.mx}}}
\end{center}

{\bf Abstract.} The operator of double differentiation on a finite interval with Robin boundary conditions perturbed by the composition of
a Volterra convolution operator and the differentia\-tion one is considered. We study the inverse problem of recovering the convolution
kernel along with a coefficient of the boundary conditions from the spectrum. We prove the uniqueness theorem and that the standard
asymptotics is a necessary and sufficient condition for an arbitrary sequence of complex numbers to be the spectrum of such an operator. A
constructive procedure for solving the inverse problem is given.



\begin{center}
{\bf 1. INTRODUCTION}
\end{center}

Consider the boundary value problem $L=L(M,h,H)$ of the form
$$
\ell y:=-y''+ \int_0^x M(x-t)y'(t)\,dt=\lambda y, \quad 0<x<\pi,       \eqno(1)
$$
$$
U(y):=y'(0)-hy(0)=0, \quad V(y):=y'(\pi)+Hy(\pi)=0,                     \eqno(2)
$$
where $\lambda$ is the spectral parameter, $M(x)$ is a complex-valued function, $(\pi-x)M(x)\in L_2(0,\pi)$ and $h,\,H\in {\mathbb C}.$

We study an inverse spectral problem for $L.$ Inverse problems of spectral analysis consist in recovering operators
from given their spectral characteristics. Such problems often appear in mathematics, mechanics, physics, electronics,
geophysics, meteorology and other branches of natural sciences and engineering. The greatest success in the inverse
problem theory has been achieved for the Sturm–Liouville operator $\ell_1y:=-y''+q(x)y$ (see, e.g., [1--4]) and
afterwards for higher-order differential operators [5--7]. For example, it is known that the potential $q(x)$ can be
uniquely determined by specifying the spectra of two boundary value problems for equation $\ell_1 y=\lambda y$ with one
common boundary condition.

For integro-differential and other classes of nonlocal operators inverse problems are more difficult for investigation,
and the classical methods either are not applicable to them or require essential modifications (see [4,\,8--19] and the
references therein). In [9] a perturbation of the Sturm-Liouville operator with Dirichlet boundary conditions by the
Volterra convolution operator was considered. It was proven that the specification of only the spectrum uniquely
determines the convolution component. Moreover, developing the idea of Borg's method a constructive procedure for
solving this inverse problem was obtained along with the local solvability and stability. In [15] the global
solvability was proved by reducing this inverse problem to solving the so-called main nonlinear integral equation,
which was solved globally. Earlier by a particular case of this approach the analogous results were obtained for the
operator (1) with Dirichlet boundary conditions [12]. Here we study the case of Robin boundary conditions. A short
version of this preprint is to appear in [20].

Let $\varphi(x,\lambda)$ be a solution of equation (1) satisfying the initial conditions
$$
\varphi(0,\lambda)=1, \quad \varphi'(0,\lambda)=h.     \eqno(3)
$$
Clearly, the eigenvalues of the problem $L$ with account of multiplicity coincide with
the zeros of the function
$$
\Delta(\lambda):=V(\varphi(x,\lambda)),     \eqno(4)
$$
which is called the characteristic function of $L.$ By the wellknown method (see, e.g., [4]) involving Rouch\'e's theorem one can prove
that the spectrum of $L$ consists of infinitely many eigenvalues $\lambda_k,$ $k\ge0.$ Moreover, the following theorem holds.

\medskip
{\bf Theorem 1. }{The spectrum $\{\lambda_k\}_{k\ge0}$ has the form
$$
\lambda_k=(k+\kappa_k)^2, \quad \{\kappa_k\}\in l_2.           \eqno(5)
$$}

As compared with the Dirichlet boundary conditions the Robin ones (2) bring additional difficulties in studying the inverse problem for
$L.$ First, we consider the following problem.

\medskip
{\bf Inverse Problem 1.} Given $\{\lambda_k\}_{k\ge0}$ and $h,\,H,$ find $M(x).$

\medskip
For this inverse problem we prove the following uniqueness theorem.

\medskip
{\bf Theorem 2. }{\it The specification of the spectrum $\{\lambda_k\}_{k\ge0}$ uniquely determines the function $M(x),$ provided that the
coefficients $h,\,H$ are known a priori.}

\medskip
We note that Inverse Problem 1 is overdetermined. Indeed, the spectrum possesses also some information on the
coefficients of the boundary conditions. In particular, we prove that if $h=0,$ then along with $M(x)$ the coefficient
$H$ is also determined, i.e. a uniqueness theorem holds for the following problem.

\medskip
{\bf Inverse Problem 2.} Given $\{\lambda_k\}_{k\ge0}$ and $h=0,$ find $H$ and
$M(x).$

\medskip
Moreover, the following theorem holds.

\medskip
{\bf Theorem 3. }{\it For arbitrary complex numbers $\lambda_k,\;k\ge0,$ of the form (5) there exists a unique (up to
values on a set of measure zero) function $M(x),$ $(\pi-x)M(x)\in L_2(0,\pi),$ and a unique number $H\in{\mathbb C},$
such that $\{\lambda_k\}_{k\ge0}$ is the spectrum of the problem $L(M,0,H).$}

\medskip
Thus, the asymptotics (5) is a necessary and sufficient condition for the solvability of Inverse Problem~2. The
importance of the assumption $h=0$ is explained in Remark~1 (see Section~4). We leave open whether the analogous
criterium can be obtained assuming only that $h$ is known a priori but $h\ne0,$ or symmetrically: $H$ is given while
$h$ is unknown.

In the next section we derive the main nonlinear integral equation of the inverse problem and prove the global
solvability of this nonlinear equation. In Section~3 we prove some auxiliary assertions along with Theorem~1. In
Section~4 we give the proof of Theorems~2 and~3, which is constructive, and provide algorithms for solving the inverse
problems (Algorithms~1 and~2).

\begin{center}
{\bf 2. MAIN NONLINEAR INTEGRAL EQUATION}
\end{center}

Let the functions $C(x,\lambda),\,S(x,\lambda)$ be solutions of equation (1) satisfying the initial conditions
$$
C(0,\lambda)=S'(0,\lambda)=1, \quad C'(0,\lambda)=S(0,\lambda)=0.
$$
Thus, according to (3), (4) we have $\varphi(x,\lambda)=C(x,\lambda)+hS(x,\lambda)$ and
$$
\Delta(\lambda)=C'(\pi,\lambda)+ hS'(\pi,\lambda)+HC(\pi,\lambda)+hHS(\pi,\lambda).      \eqno(6)
$$

Let $\rho^2=\lambda.$ The following representation is wellknown:
$$
S(x,\lambda)=\frac{\sin\rho x}{\rho}+\int_0^x P(x,t)\frac{\sin\rho(x-t)}{\rho}\,dt, \quad 0\le x\le\pi, \eqno(7)
$$
where $P(x,t)$ is the kernel of the transformation operator. In [12] it was shown that
$$
P(x,t)=\sum_{\nu=1}^\infty \frac{(x-t)^\nu}{\nu!} N^{*\nu}(t),  \eqno(8)
$$
where
$$
N^{*1}(x)=N(x), \quad N^{*(\nu+1)}(x)=N*N^{*\nu}(x)=\int_0^x N(x-t)N^{*\nu}(t)\,dt,
\;\; \nu\ge1,
$$
and the function $N(x),\,(\pi-x)N(x)\in L_2(0,\pi),$ is connected with $M(x)$ by the relation
$$
M(x)=2N(x)-\int_0^x dt \int_0^t N(t-\tau)N(\tau)\,d\tau, \quad 0<x<\pi. \eqno(9)
$$

The following lemma gives further representations, which we use in the sequel.

\medskip
{\bf Lemma 1.} {\it The following representations hold:
$$
S(x,\lambda)=K(x,x)-\rho\int_0^x K(x,t)\sin\rho(x-t)\,dt,  \eqno(10)
$$
$$
S'(x,\lambda)=R(x,x)-\rho\int_0^x R(x,t)\sin\rho(x-t)\,dt, \eqno(11)
$$
$$
C(x,\lambda)=1-\rho\int_0^x Q(x,t)\sin\rho(x-t)\,dt, \eqno(12)
$$
$$
C'(x,\lambda)=-\rho\sin\rho x-\rho\int_0^x P(x,t)\sin\rho(x-t)\,dt, \eqno(13)
$$
where
$$
K(x,t)=t+\int_0^t (t-\tau)P(x,\tau)\,d\tau,     \eqno(14)
$$
$$
R(x,t)=1+\int_0^t P(x,\tau)\,d\tau+\int_0^t (t-\tau)P_x(x,\tau)\,d\tau, \eqno(15)
$$
$$
Q(x,t)=1+\int_0^t P(x-t+\tau,\tau)\,d\tau.     \eqno(16)
$$}

\medskip
{\it Proof.} The integration by parts in (7) gives
$$
S(x,\lambda)=\int_0^x \Big(1+\int_0^t P(x,\tau)\,d\tau\Big)\cos\rho(x-t)\,dt,      \eqno(17)
$$
Integrating (17) by parts we arrive at (10) and (14). Then differentiating (10) with respect to $x$ we get
$$
S'(x,\lambda)=\frac{d}{dx}K(x,x)-\rho\int_0^x K_x(x,t)\sin\rho(x-t)\,dt-\lambda\int_0^x K(x,t)\cos\rho(x-t)\,dt.
$$
Integrating by parts in the last term we arrive at
$$
S'(x,\lambda)=\frac{d}{dx}K(x,x)-\rho\int_0^x (K_x(x,t)+K_t(x,t))\sin\rho(x-t)\,dt.
$$
Taking (14) and (15) into account we get (11).

Further, we have
$$
C(x,\lambda)=1-\lambda\int_0^x S(t,\lambda)\,dt.         \eqno(18)
$$
Indeed, put
$$
v(x):=S'(x,\lambda)-\int_0^x M(x-t)S(t,\lambda)\,dt, \quad u(x):=1-\lambda\int_0^x S(t,\lambda)\,dt.
$$
Since $v(0)=u(0)$ and $v'(x)=u'(x),$ we have $v(x)=u(x).$ On the other hand, since
$$
u'(x)=-\lambda S(x,\lambda), \quad u''(x)=-\lambda S'(x,\lambda),
$$
we get
$$
\ell u=\lambda\Big(S'(x,\lambda)-\int_0^x M(x-t)S(t,\lambda)\,dt\Big)=\lambda v(x)=\lambda u(x).
$$
Taking into account $u(0)=C(0,\lambda)$ and $u'(0)=C'(0,\lambda)$ we arrive at (18).

Differentiating (18) and substituting (7) therein we get (13).

Finally, substituting (7) into (18) we get
$$
C(x,\lambda)=1-\rho\int_0^x\sin\rho t\,dt-\rho\int_0^x dt\int_0^t P(t,t-\tau)\sin\rho\tau\,d\tau.
$$
Changing the order of integration yields
$$
C(x,\lambda)=1-\rho\int_0^x\Big(1+\int_{x-t}^x P(\tau,\tau-x+t)\,d\tau\Big)\sin\rho(x-t)\,dt,
$$
which gives (12) and (16). $\hfill\Box$

\medskip
The next lemma is a direct corollary of formulae (6), (10)--(13).

\medskip
{\bf Lemma 2. }{\it The characteristic function has the form
$$
\Delta(\lambda)=-\rho\sin\rho\pi+\alpha+\rho\int_0^\pi w(x)\sin\rho x\,dx, \quad w(x)\in L_2(0,\pi). \eqno(19)
$$
Here
$$
\alpha=hR(\pi,\pi;N)+H+hHK(\pi,\pi;N), \eqno(20)
$$
$$
-w(\pi-x)=P(\pi,x;N)+hR(\pi,x;N)+HQ(\pi,x;N)+hHK(\pi,x;N), \eqno(21)
$$
where we add the argument ''$N$'' in order to indicate the dependence on $N(x).$}

\medskip
The relation (21) can be considered as a nonlinear equation with respect to $N(x).$ We call it {\it main nonlinear
integral equation} of the inverse problem. The main equation (21) can be rewritten in the explicit form. Indeed, by
virtue of (8), (14)--(16) we have
$$
P(\pi,x;M)=\sum_{\nu=1}^\infty \frac{(\pi-x)^\nu}{\nu!} N^{*\nu}(x),  \eqno(22)
$$
$$
R(\pi,x;M)=1+\sum_{\nu=1}^\infty \frac{1}{\nu!}\int_0^x (\pi-t)^\nu N^{*\nu}(t)\,dt +\sum_{\nu=1}^\infty
\frac{1}{\nu!}\int_0^x \nu(x-t)(\pi-t)^{\nu-1} N^{*\nu}(t)\,dt,  \eqno(23)
$$
$$
Q(\pi,x;M)=1+\sum_{\nu=1}^\infty \frac{(\pi-x)^\nu}{\nu!}\int_0^x  N^{*\nu}(t)\,dt,  \eqno(24)
$$
$$
K(\pi,x;M)=x+\sum_{\nu=1}^\infty \frac{1}{\nu!}\int_0^x(x-t)(\pi-t)^\nu N^{*\nu}(t)\,dt.  \eqno(25)
$$
Substituting (22)--(25) into (21) we get
$$
f(x)=\sum_{\nu=1}^\infty\Big(\psi_\nu(x) N^{*\nu}(x) +\int_0^x\Psi_\nu(x,t) N^{*\nu}(t)\,dt\Big), \eqno(26)
$$
where $f(x)=-w(\pi-x)-h-H-hHx,$ $\psi_\nu(x)=(\pi-x)^\nu/\nu!$ and
$$
\Psi_\nu(x,t)=\frac{1}{\nu!}\Big(H(\pi-x)^\nu+h(\pi-t)^{\nu-1}\Big(\pi-t+(x-t)(\nu+H(\pi-t))\Big)\Big).
$$

\medskip
{\bf Theorem 4. }{\it For each function $f(x)\in L_2(0,\pi)$ and any complex numbers $h,\,H$ equation (26) has a unique
solution $N(x),\,(\pi-x)N(x)\in L_2(0,\pi).$}

\medskip
{\it Proof.} By virtue of Theorem 4 in [10], equation (26) has a unique solution $N(x),$ which belongs to $L_2(0,T)$
for each $T\in(0,\pi).$ Following [10] we represent this solution in the form
$$
N(x)=N_1(x)+N_2(x),
$$
where $N_1(x)\in L_2(0,\pi)$ and $N_2(x)=0$ on $(0,\pi/2).$ Then
we have
$$
N^{*\nu}(x)=N_1^{*\nu}(x)+\nu N_1^{*(\nu-1)}*N_2(x), \quad \nu\ge2.
$$
Substituting this into (26) we arrive at
$$
f(x)-\mu_1(x) = (\pi-x)N_2(x) + \int_{\pi/2}^x A(x,t)\,(\pi-t)N_2(t)\,dt, \quad \frac\pi2<x<\pi,
$$
where
$$
\mu_1(x)=\sum_{\nu=1}^\infty\Big(\psi_\nu(x) N_1^{*\nu}(x) +\int_0^x\Psi_\nu(x,t)
N_1^{*\nu}(t)\,dt\Big),
$$
$$
A(x,t)=\frac{1}{\pi-t}\Big(\Psi_1(x,t)+\sum_{\nu=2}^\infty\nu \Big(\psi_\nu(x)
N_1^{*(\nu-1)}(x-t) +\int_t^x\Psi_\nu(x,\tau)
N_1^{*(\nu-1)}(\tau-t)\,d\tau\Big)\Big)
$$
are square-integrable functions. Hence $(\pi-x)N_2(x)\in L_2(0,\pi).$ $\hfill\Box$

\newpage

\begin{center}
{\bf 3. AUXILIARY ASSERTIONS}
\end{center}

In this section for convenience of the reader we prove Theorem~1 and further auxiliary assertions, which we use in
Section~4 for solving the inverse problems. We note that Theorem~1 is a direct corollary of the following assertion.

\medskip
{\bf Lemma 3. } {\it Any function of the form (19) has infinitely many zeros $\lambda_k,\,k\ge0,$ having the
asymptotics (5).}

\medskip
{\it Proof.} For $\rho\in G_\delta:=\{\rho:|\rho-k|\ge \delta,\,k\in{\mathbb Z}\},$ $\delta>0,$ we have
$|\sin\rho\pi|\ge C_\delta\exp(|{\rm Im}\,\rho|\pi),$ $C_\delta>0.$ Therefore, for sufficiently large $|\rho|,$
$\rho\in G_\delta,$ the following inequality holds
$$
|\rho\sin\rho\pi|>\Big|\alpha+\rho\int_0^\pi w(x)\sin\rho x\,dx\Big|.
$$
According to Rouch\'e's theorem there are exactly $N+1$ zeros $\lambda_k,\;k=\overline{0,N},$ of the function
$\Delta(\lambda)$ lying inside the contour $\Gamma_N=\{\lambda:\; |\lambda|=(N+1/2)^2\}$ for sufficiently large $N.$
Moreover, for each $\delta>0$ there exists $k_\delta$ such that for $|k|>k_\delta$ there is exactly one zero $\rho_k$
of the function $\Delta(\rho^2)$ inside the contour $\gamma_k(\delta)=\{\rho:\;|\rho-k|=\delta\}.$ Thus,
$\rho_k=\sqrt{\lambda_k}= k+\kappa_k,$ where $\kappa_k=o(1).$ Substituting this into (19) we get $\{\kappa_k\}\in l_2$
and (5) is proved. $\hfill\Box$

\medskip
Analogously to Theorem~1.1.4 in [4] using Hadamard's factorization theorem one can obtain the following assertion.

\medskip
{\bf Lemma 4. }{\it The function $\Delta(\lambda)$ is uniquely determined by its zeros by the formula
$$
\Delta(\lambda)=\pi(\lambda_0-\lambda)\prod_{k=1}^\infty\frac{\lambda_k-\lambda}{k^2}. \eqno(27)
$$}

{\it Proof. } It follows from (19) that $\Delta(\lambda)$ is entire in $\lambda$ of order 1/2, and consequently by
Hadamard`s factorization  theorem, $\Delta(\lambda)$ is uniquely determined up to a multiplicative constant by its
zeros:
$$
\Delta(\lambda) = C\lambda^s\prod_{\lambda_k\ne0}\Big(1-\frac{\lambda}{\lambda_k}\Big),          \eqno(28)
$$
where $s\ge0$ is the multiplicity of the eigenvalue $\lambda=0.$ Consider the function
$$
\Delta_0(\lambda) := -\rho\sin\rho\pi =-\lambda\pi \prod_{k=1}^{\infty} \Big(1-\frac{\lambda}{k^2}\Big).
$$
Then
$$
\frac{\Delta(\lambda)}{\Delta_0(\lambda)} = C\frac{\lambda^s}{\pi} \prod_{\lambda_k=0} \frac{n_k}{k^2-\lambda}
\prod_{\lambda_k\ne0} \frac{n_k}{\lambda_k} \prod_{\lambda_k\ne0}\Big(1+\frac{\lambda_k-k^2}{k^2-\lambda}\Big), \quad
{\rm where} \quad n_k=\left\{\begin{array}{l}k^2,\;k\in{\mathbb N},\\[1mm]1,\;k=0.\end{array}\right.
$$
On the other hand, taking (5) and (19) into account we calculate
$$
\lim_{\lambda \to -\infty} \frac{\Delta(\lambda)}{\Delta_0(\lambda)} =1, \quad \lim_{\lambda \to -\infty}
\prod_k\Big(1+\frac{\lambda_k-k^2} {k^2-\lambda}\Big)=1,
$$
and hence
$$
C=(-1)^s\pi\prod_{\lambda_k=0} \frac{1}{n_k}\prod_{\lambda_k\ne0} \frac{\lambda_k}{n_k}.
$$
Substituting this into (28) we arrive at (27). $\hfill\Box$

\medskip
Thus, according to Lemmas 3 and 4 any function of the form (19) is uniquely determined by formula (27) from its zeros,
which, in turn, have the asymptotics (5). By the standard approach (see, e.g., Lemma 3.3 in [12]) one can prove the
following inverse assertion.

\medskip
{\bf Lemma 5. }{\it (i) Let arbitrary complex numbers $\lambda_k,\;k\ge0,$ of the form (5) be given. Then the function
$\Delta(\lambda)$ determined by (27) has the form (19) with certain $\alpha\in{\mathbb C}$ and complex-valued function
$w(x)\in L_2(0,\pi).$}

\medskip
{\it Proof.} Once again using the function
$$
\Delta_0(\lambda):=-\rho\sin\rho\pi=-\pi\lambda\prod_{k=1}^{\infty} \frac{k^2-\lambda}{k^2}, \eqno(29)
$$
from (27) we get
$$
\Delta(\lambda)=\Delta_0(\lambda)F(\lambda), \quad F(\lambda) = \frac{\lambda-\lambda_0}{\lambda}\prod_{k=1}^{\infty}
\frac{\lambda_k-\lambda}{k^2-\lambda}.  \eqno(30)
$$
Put
$$
\beta_k:=\Delta(k^2), \quad \theta_k:=\frac{\beta_k-\alpha}{k}, \quad k\ge1. \eqno(31)
$$
By virtue of (29), (30), we have
$$
\beta_k=(-1)^{k+1}\frac{\pi}{2k^2}(\lambda_0-k^2)(\lambda_k-k^2)d_k, \quad d_k=\prod_{\nu=1,\;\nu\ne k}^\infty
\frac{\lambda_\nu-k^2}{\nu^2-k^2}, \quad k\ge1. \eqno(32)
$$
Hence, according to (5), (31), (32) we get
$$
\theta_k=(-1)^{k+1}\frac{\pi}{2k^3}(\lambda_0-k^2)(2\kappa_k k+\kappa_k^2)d_k-\frac{\alpha}{k}.
$$
Since the sequence $\{d_k\}$ is bounded, we have $\{\theta_k\}\in l_2.$ Determine the function $w(x)\in L_2(0,\pi)$
such that
$$
\theta_k=\int_0^\pi w(x)\sin kx\,dx, \quad k\ge1,  \eqno(33)
$$
and denote
$$
\theta(\rho)=\int_0^\pi w(x)\sin\rho x\,dx.
$$
Consider the function
$$
S(\rho):=\frac{\theta(\rho)-\frac{\Delta(\rho^2)-\alpha}{\rho}-\sin\rho\pi}{\sin\rho\pi}
=\frac{\theta(\rho)}{\sin\rho\pi}+F(\rho^2)+\frac{\alpha}{\rho\sin\rho\pi}-1,
$$
which is, by (31) and (33) after removing the singularities, entire in $\rho.$ Let us show that $|F(\rho^2)|<C_\delta$
for $\rho\in G_\delta=\{\rho:\;|\rho-k|\ge \delta, \, k\in{\mathbb Z}\},$ $\delta>0.$ Denote $\kappa_{-k}:=-\kappa_k,\;
k\in{\mathbb N}.$ Then
$$
F(\rho^2)=\frac{\rho^2-\lambda_0}{\rho^2}\prod_{|k|\in{\mathbb N}}\Big(1+\frac{\kappa_k}{k-\rho}\Big).
$$
For a fixed $\delta>0$ choose $N$ such that $|\kappa_k|\le\delta/2$ for $k\ge N.$ We have for $\rho\in G_\delta$
$$
F(\rho^2)=\exp(H_N(\rho))\frac{\rho^2-\lambda_0}{\rho^2}\prod_{0<|k|<N}\Big(1+\frac{\kappa_k}{k-\rho}\Big), \eqno(34)
$$
where
$$
H_N(\rho)=\sum_{|k|\ge N} \ln \Big(1+\frac{\kappa_k}{k-\rho}\Big)= \sum_{|k|\ge N}
\frac{\kappa_k}{k-\rho}\,\sum_{\nu=0}^\infty \frac{(-1)^\nu}{\nu +1} \Big(\frac{\kappa_k}{k-\rho}\Big)^\nu.
$$
For $\rho\in G_\delta$ using the Cauchy-Bunyakovsky inequality we arrive at
$$
|H_N(\rho)|\le  C\Big(\sum_{|k|\ge N} \frac{1}{|k-\rho|^2}\Big)^{1/2}
$$
and consequently $H_N(\rho)$ is bounded in $G_\delta.$ Thus, from (34) it follows that $|F(\rho^2)|< C_\delta$ for
$\rho\in G_\delta.$ Moreover, it is easy to see that $F(\rho^2)\to 1$ as $|Im\rho|\to\infty.$ Using the maximum modulus
principle we conclude that $S(\rho)$ is bounded and consequently, according to Liouville's theorem, it is constant.
Since $S(\rho)\to 0$ for $|Im\rho|\to\infty,$ $S(\rho)\equiv 0$ and formula (19) holds. $\hfill\Box$

\begin{center}
{\bf 4. SOLUTION OF THE INVERSE PROBLEM}
\end{center}

By virtue of (19), (27), the number $\alpha$ can be found by the formula
$$
\alpha=\pi\lambda_0\prod_{k=1}^\infty\frac{\lambda_k}{k^2}.     \eqno(35)
$$
According to (19) the function $w(x)$ can be reconstructed as the Fourier series
$$
w(x)= \frac2\pi\sum_{k=1}^\infty\frac{\Delta(k)}{k}\sin kx - \alpha\frac{\pi-x}{\pi}. \eqno(36)
$$

\medskip
{\it Proof of Theorem 2.} According to Lemmas 2 and 4 the specification of the spectrum uniquely determines the
function $w(x).$ Since $h$ and $H$ are given, by virtue of Theorem~4 the function $N(x)$ is a unique solution of the
main equation (26). Hence, the function $M(x)$ is determined uniquely and can be reconstructed by formula (9).
$\hfill\Box$

\medskip
Let the spectrum $\{\lambda_k\}_{k\ge0}$ of a problem $L(M,h,H)$ along with the
numbers $h,\,H$ be given. Then the function $M(x)$ can be constructed by the
following algorithm.

\medskip
{\bf Algorithm 1.} (i) Construct the function $w(x)$ by formulae (27), (35) and (36);

(ii) find $N(x)$ by solving the main equation (26);

(iii) construct $M(x)$ by formula (9).

\medskip
{\it Proof of Theorem 3.} Using the given numbers $\lambda_k$ we construct the function $\Delta(\lambda)$ by formula
(27). According to Lemma~5 it has the representation (19) with certain number $\alpha$ and function $w(x)\in
L_2(0,\pi).$ Let $N(x)$ be the solution of the main equation (26) with this function $w(x)$ and $h=0,\,H=\alpha.$
Determine the function $M(x),\,(\pi-x)M(x)\in L_2(0,\pi),$ by formula (9) and consider the corresponding boundary value
problem $L=L(M,0,H).$ It is easy to see that the constructed function $\Delta(\lambda)$ is the characteristic function
of this problem~$L.$ Thus, the spectrum of the latter coincides with $\{\lambda_k\}_{k\ge0}.$ According to (9), (19)
and~(27) the uniqueness of $M(x)$ follows from the uniqueness of the solution of equation (26). $\hfill\Box$

\medskip
Let a sequence $\{\lambda_k\}_{k\ge0}$ of the form (5) be given. According to Theorem~3 there exists a unique boundary
value problem $L(M,0,H)$ with the spectrum $\{\lambda_k\}_{k\ge0},$ which can be constructed by the following
algorithm.

\medskip
{\bf Algorithm 2.} (i) Having calculated the number $\alpha$ by formula (35), put
$$
h:=0, \quad H:=\alpha;
$$
and construct the function $w(x)$ by formulae (27), (36);

(ii) find $N(x)$ by solving the main equation (26);

(iii) construct $M(x)$ by formula (9).

\medskip
{\bf Remark 1.} The importance of the assumption $h=0$ can be seen from formula (20). Indeed, since $N(x)$ is unknown,
the number $\alpha$ being, in turn, uniquely determined by the spectrum determines $H$ only if $h=0.$

\medskip
{\bf Acknowledgments.} The first author was supported in part by Russian Foundation
for Basic Research (Grants 15-01-04864, 13-01-00134) and by the Ministry of
Education and Science of RF (Grant 1.1436.2014K). The second author was supported by
CONACYT grant No. 153184, CIC–UMSNH, M\'exico

\begin{center}
{\bf REFERENCES}
\end{center}

\begin{enumerate}
\item[{[1]}] Borg G. {\it Eine Umkehrung der Sturm-Liouvilleschen Eigenwertaufgabe}, Acta Math. 78 (1946), 1--96.

\item[{[2]}] Marchenko V.A. {\it Sturm-Liouville Operators and Their Applications}, Naukova Dum\-ka, Kiev, 1977;
English transl., Birkh\"auser, 1986.

\item[{[3]}] Levitan B.M. {\it Inverse Sturm-Liouville Problems}, Nauka, Moscow, 1984; English transl., VNU
Sci.Press, Utrecht, 1987.

\item[{[4]}] Freiling G. and Yurko V.A. {\it Inverse Sturm-Liouville Problems and Their Applications}, NOVA Science
Publishers, New York, 2001.

\item[{[5]}] Beals R., Deift P., Tomei C. {\it Direct and Inverse Scattering on the Line}, Mathematica Surveys and
Monographs, 28. AMS, Providence, RI, 1988.

\item[{[6]}] Yurko V.A. {\it Inverse Spectral Problems for Differential Operators and Their Applica\-tions}, Gordon
and Breach Science Publishers, Amsterdam, 2000.

\item[{[7]}] Yurko V.A. {\it Method of Spectral Mappings in the Inverse Problem Theory}, Inverse and Ill-posed
Problems Series. VSP, Utrecht, 2002.

\item[{[8]}] Eremin M.S. {\it An inverse problem for a second-order integro-differential equation with a singularity}, Diff.
Uravn. 24 (1988), no.2, 350--351.

\item[{[9]}] Yurko V.A. {\it An inverse problem for integro-differential operators}, Mat. Zametki, 50 (1991), no.5,
134--146 (Russian); English transl. in Math. Notes 50 (1991), no. 5--6, 1188--1197.

\item[{[10]}] Buterin S.A. {\it The inverse problem of recovering the Volterra convolution operator from the incomplete
spectrum of its rank-one perturbation}, Inverse Problems 22 (2006), 2223--2236.

\item[{[11]}] Buterin S.A. {\it Inverse spectral reconstruction problem for the convolution operator pertur\-bed by a one-dimensional
operator}, Mat. Zametki 80 (2006), no.5, 668--682 (Russian); English transl. in Math. Notes 50 (2006), no. 5, 631--644.

\item[{[12]}] Buterin S.A. {\it On an inverse spectral problem for a convolution integro-differential operator}, Res. Math. 50 (2007),
no.3-4, 73--181.

\item[{[13]}] Albeverio S., Hryniv R.O. and Nizhnik L.P. {\it Inverse spectral problems for non-local Sturm-\-Liouville operators}, Inverse
Problems 23 (2007), 523--535.

\item[{[14]}] Kuryshova Ju.V. {\it Inverse spectral problem for integro-differential operators}, Mat. Zametki 81 (2007),
no.6, 855--866; English transl. in Math. Notes 81 (2007), no.6, 767--777.

\item[{[15]}] Buterin S.A. {\it On the reconstruction of a convolution perturbation of the Sturm-Liouville operator from the spectrum},
Diff. Uravn. 46 (2010), 146–149 (Russian); English transl. in Diff. Eqns. 46 (2010), 150–-154.

\item[{[16]}] Kuryshova Yu.V. and Shieh C.-T. {\it An inverse nodal problem for integro-differential opera\-tors}, J. Inverse Ill-posed
Probl. 18 (2010), no.4, 357--369.

\item[{[17]}] Nizhnik L.P. {\it Inverse spectral nonlocal problem for the first order ordinary differential equation}, Tamkang J. Math.
42 (2011), no.3, 385--394.

\item[{[18]}] Yang C.-F. and Yurko V.A. {\it Recovering differential operators with nonlocal boundary conditions}, arXiv:1410.2017v1 [math.SP],
2014, 12pp.

\item[{[19]}] Yurko V.A. {\it An inverse spectral problems for integro-differential operators}, Far East J. Math. Sci. 92 (2014), no.2,
247--261.

\item[{[20]}] Buterin S.A. and Choque Rivero A.E. {\it On inverse problem for a convolution integro-differential operator with Robin
boundary conditions}, Applied Math. Letters (to appear)

\end{enumerate}

\end{document}